\documentclass{amsart}

\newtheorem{theorem}{Theorem}[section]
\newtheorem*{maintheorem}{Theorem}
\newtheorem{lemma}[theorem]{Lemma}
\newtheorem{proposition}[theorem]{Proposition}
\newtheorem{corollary}[theorem]{Corollary}

\theoremstyle{definition}

\newtheorem{remark}[theorem]{Remark}
\newtheorem*{acknowledgement}{Acknowledgement}

\theoremstyle{remark}

\usepackage{amscd,amssymb}
\newcommand\mynote[1]{\marginpar{\ \\ \small \tt #1}}
\newcommand\bel[1]{{\mynote{#1}}\begin{equation}\label{#1}}
\newcommand\mylabel[1]{\label{#1}}

\newcommand{\ZZ}{\mathbb{Z}}
\newcommand{\QQ}{\mathbb{Q}}

\newcommand{\CC}{\mathbb{C}}

\newcommand{\PP}{\mathbb{P}}

\newcommand  {\shA}     {\mathcal{A}}

\newcommand  {\shB}     {\mathcal{B}}
\newcommand  {\shC}     {\mathcal{C}}

\newcommand  {\shEnd}   {\mathcal{E}\!\text{\textit{nd}}}

\newcommand  {\shE}     {\mathcal{E}}
\newcommand  {\shF}     {\mathcal{F}}
\newcommand  {\shG}     {\mathcal{G}}
\newcommand  {\shH}     {\mathcal{H}}

\newcommand  {\shL}     {\mathcal{L}}

\newcommand  {\foH}     {\mathfrak{H}}

\newcommand  {\foU}     {\mathfrak{U}}

\newcommand  {\ab}      {{\operatorname{ab}}}

\newcommand  {\an}      {{\text{an}}}

\newcommand  {\Aut}     {\operatorname{Aut}}

\newcommand  {\Br}      {\operatorname{Br}}

\newcommand  {\cl}      {\operatorname{cl}}
\newcommand  {\cH}      {\check{H}}
\newcommand  {\coker}   {\operatorname{coker}}

\newcommand  {\GL}      {\operatorname{GL}}

\newcommand  {\Hom}     {\operatorname{Hom}}

\newcommand  {\id}      {\operatorname{id}}
\newcommand  {\im}      {\operatorname{im}}

\newcommand  {\dirlim}  {\varinjlim}
\newcommand  {\invlim}  {\varprojlim}

\newcommand  {\lra}     {\longrightarrow}

\newcommand  {\Mat}     {\operatorname{Mat}}

\renewcommand{\O}       {\mathcal{O}}

\newcommand  {\op}      {{\operatorname{op}}}

\newcommand  {\Pic}     {\operatorname{Pic}}
\newcommand  {\PGL}     {\operatorname{PGL}}

\newcommand  {\ra}      {\rightarrow}

\newcommand  {\SL}      {\operatorname{SL}}

\newcommand  {\Tor}     {\operatorname{Tor}}

\def\mydate{\number\day\space\ifcase\month \or January\or February\or March\or
April\or May\or June\or July\or
August\or September\or October\or November\or December\fi \space\number\year}
\begin{document}

\title[Complex-analytic surfaces]
      {Topological methods  for complex-analytic Brauer groups}

\author[Stefan Schroer]{Stefan Schr\"oer}
\address{Mathematisches Institut, Heinrich-Heine-Universit\"at,
40225 D\"usseldorf, Germany}
\curraddr{}
\email{schroeer@math.uni-duesseldorf.de}


\dedicatory{Final version, 7 January 2005}

\begin{abstract}
Using methods from algebraic topology and
group cohomology, I pursue Grothendieck's question on
equality of geometric and  cohomological Brauer
groups in the context of complex-analytic spaces.
The main result is that equality holds under
suitable assumptions on the fundamental group
and the Pontrjagin dual of the second homotopy group.
I apply this to Lie groups, Hopf manifolds,
and complex-analytic surfaces.
\end{abstract}

\maketitle
\tableofcontents

%===========================================================
\section*{Introduction}
\mylabel{introduction}

The goal of this paper is to pursue Grothendieck's question
on Brauer groups in the context of
complex-analytic spaces, using methods from
algebraic topology and group cohomology.
This yields new results on Lie groups, Hopf manifolds,
and surfaces.

Let me  recall  Grothendieck's question. Suppose $X$ is a
topological space endowed with a sheaf of rings $\O_X$. The
\emph{cohomological Brauer group} $\Br'(X)$ is defined as the torsion
part of the cohomology group $H^2(X,\O_X^\times)$. One likes to
have a geometric interpretation of such cohomology classes. A
possible interpretation is in terms of principal
$\PGL_r$-bundles $P\ra X$, via the nonabelian coboundary map
$H^1(X,\PGL_r(\CC))\ra H^2(X,\O_X^\times)$. The group of
equivalence classes of principal $\PGL_r$-bundles is called the
\emph{Brauer group} $\Br(X)$. The coboundary map yields a canonical inclusion
$\Br(X)\subset\Br'(X)$, and Grothendieck \cite{GB II} asked whether
this inclusion is actually a bijection. This is a major open
problem in the theory of Brauer groups.

In algebraic topology, Serre solved Grothendieck's question if $X$ is
a finite CW-complex and the sheaf of rings is the sheaf of continuous complex
functions $\shC_X$. It turns out that here  $\Br(X)=\Br'(X)$  equals  the torsion
part of $H^3(X,\ZZ)$. In contrast, there are only few general
results in algebraic geometry. In this context
$X$ is an algebraic scheme and $\O_X$ is its structure sheaf. 
Gabber's Theorem   tells us that
$\Br(X)=\Br'(X)$ for any affine scheme \cite{Gabber 1981}.  Grothendieck himself
showed that equality holds for smooth algebraic surfaces
\cite{GB II}, and  I treated the case of algebraic surfaces with isolated
singularities \cite{Schroeer 2001}.

Grothendieck's question shows up in various areas.
To mention a few: In moduli theory, Brauer groups are used
in order to determine whether a coarse moduli space is actually a fine
moduli space. In stack theory, Brauer groups are important to
detect quotient stacks, as explained in the work of
Edidin, Hassett, Kresch, and Vistoli \cite{Edidin; Hassett; Kresch; Vistoli 1999}.
In Homological Mirror Symmetry, Brauer groups are used for twisting
derived categories.

This paper deals with complex-analytic spaces, which are not
necessarily algebraic. Here methods of algebraic geometry
frequently  break down, largely due to extension problems
involving coherent sheaves. In this context there are  two general
results: Elencwajg and  Narasimhan  showed $\Br(X)=\Br'(X)$ for
complex tori \cite{Elencwajg; Narasimhan 1983}, and Huybrechts and
myself recently proved it for complex-analytic K3-surfaces
 with methods from
differential geometry \cite{Huybrechts; Schroeer 2003}.
Here I prove a general result on complex-analytic Brauer groups that  depends only
on the homotopy type of the underlying topological space.

\begin{maintheorem}
Let $X$ be a complex-analytic space.
Suppose $\pi_1(X)$ is a good group in Serre's sense,
and that the  subgroup of $\pi_1(X)$-invariants inside
the Pontrjagin dual $\Hom(\pi_2(X),\QQ/\ZZ)$
is trivial. Then the inclusion $\Br(X)\subset\Br'(X)$ is an equality.
\end{maintheorem}

Serre introduced the notion of \emph{good groups},
which has to do with profinite completions,
in the context of Galois cohomology \cite{Serre 1972}.
I will recall this somewhat technical concept in Section \ref{serre's good groups}.
Note that free groups and  polycyclic groups are good.

The conditions in the theorem appear bizarre, but it applies
directly to complex Lie groups and Hopf manifolds:

\begin{maintheorem}
Let $X$ be a complex Lie group or a Hopf manifold.
Then the inclusion $\Br(X)\subset\Br'(X)$ is an equality.
\end{maintheorem}

This generalizes results of
Iversen on characterfree algebraic groups \cite{Iversen 1976},
and Hoobler \cite{Hoobler 1972}, Berkovi\v{c} \cite{Berkovic 1972}, and
Elencwajg and Narasimhan \cite{Elencwajg; Narasimhan 1983} on abelian
varieties and complex tori.
Turning to surfaces,  we obtain the second main result of this paper:

\begin{maintheorem}
\mylabel{main result}
Let $S$ be a smooth compact complex-analytic surface
with $b_1\neq 1$.
Then the inclusion $\Br(S)\subset \Br'(S)$ is an equality.
\end{maintheorem}

Here the main challenge is the case of elliptic surfaces.
Usually, such surfaces do not satisfy the required conditions
on $\pi_1(S)$ and $\pi_2(S)$, due to the presence of singular fibers
in the elliptic fibration $S\ra B$.
However, there is always a Zariski open subset $U\subset S$
with the desired properties. Some additional arguments
then show that this is enough for our purpose.
A key ingredient is Hoobler's result \cite{Hoobler 1972}, see also Gabber
\cite{Gabber 1981}, that any cohomology class $\beta\in\Br'(X)$
mapping into  $\Br(Y)$ for some finite flat covering $Y\ra X$ lies
in the Brauer group.

My results for surfaces with $b_1=1$ are less definite.
Such surfaces are also called of class VII.
To date, this is the only class of surface resisting complete
classification.
Any known surface of class VII blows down to one of
the three following types: Hopf surfaces, Inoue surfaces,
and surfaces containing a global spherical shell.
The latter is a holomorphic embedding of a thickened
3-sphere with connected complement.

\begin{maintheorem}
We have $\Br(S)=\Br'(S)$ for any surface with $b_1=1$ whose minimal model
is either a  Hopf surface, Inoue surface,
or a surface containing a global spherical shell.
\end{maintheorem}

According to the \emph{GSS-conjecture}, any surface of class VII
should belong to one of these classes. If true, this
would complete the Kodaira classification.
The work of Dloussky, Oeljeklaus, and Toma gives considerable
positive evidence
\cite{Dloussky; Oeljeklaus; Toma 2000}, \cite{Dloussky; Oeljeklaus; Toma 2003}.
On the other hand, there are almost no results on
fundamental groups of hypothetical surfaces of class VII.
A notable exception is the work of Carlson and
Toledo on representations of class VII fundamental group
in fundamental groups of hyperbolic Riemannian manifolds
\cite{Carlson; Toledo 1997}.

Here is a plan for the paper:
In Section \ref{analytic brauer}, I set down some definitions
concerning analytic Brauer groups, and also establish some
Gaga type facts.
In Section \ref{nonpurity results} we shall see that
analytic Brauer groups might differ strongly from algebraic
Brauer groups on noncompact surfaces. I relate this to
the Shafarevich Conjecture, and give an application
regarding  the existence of nonalgebraic vector bundles on
pointed algebraic surfaces.
Section \ref{serre's good groups} contains a discussion of
Serre's notion of good groups.
I use good groups to prove $\Br(X)=\Br'(X)$ for certain complex-analytic spaces in
Section \ref{applications to complex spaces}. As application I discuss
the case of complex Lie groups and Hopf manifolds.
To apply the results to
complex-analytic surfaces, we still have to improve them.
This happens in Section \ref{from open to compact surfaces}.
In Section \ref{elliptic surfaces}, we then solve the
case of elliptic surfaces.
In Section \ref{Surfaces of class VII}, I analyze the
case of surfaces of class VII.
In Section \ref{main results}, I combine our result with
already known results. There is also a  discussion
of open problems.

\begin{acknowledgement}
I thank Ingrid Bauer, Fr\'ed\'eric Campana, Fabrizio Catanese, Raymond Hoobler,  Uwe Jannsen, and
Thomas Peternell  for stimulating discussions.
I also thank the referee, who pointed out some mistakes, and whose comments
helped to improve the exposition.
Moreover, I thank Le Van Ly for careful proofreading.
\end{acknowledgement}

%===========================================================
\section{Analytic Brauer groups}
\mylabel{analytic brauer}

In this section we shall introduce  notation and establish
some useful facts  on analytic Brauer groups. Throughout, $X$ denotes a
complex-analytic space, and $\O_X$ is its  sheaf of holomorphic
functions.
One way to define Brauer groups is in terms of
holomorphic principal $\PGL_r(\CC)$-bundles, which is well-suited
for our purposes. Examples of such bundles are projectivisations
$P=\PP(\shE)$ of locally free $\O_X$-modules $\shE$ of rank $r>0$.
The Brauer group measures to which extent there are other
principal bundles as follows:

Suppose $P\ra X$ is a principal $\PGL_r(\CC)$-bundle, and $P'\ra
X$ is a  principal $\PGL_{r'}(\CC)$-bundle. Using the homomorphism
$$ \PGL_r(\CC)\times\PGL_{r'}(\CC)\lra\PGL_{rr'}(\CC),\quad
(A,A')\longmapsto A\otimes A' $$
we obtain another principal
$\PGL_{rr'}(\CC)$-bundle $P\otimes P'$. One says that $P$ and
$P'$ are equivalent if there are locally free $\O_X$-modules
$\shE',\shE$ of rank $r',r>0$ so that $P\otimes\PP(\shE')$ and
$P'\otimes\PP(\shE)$ are isomorphic. The group of equivalence
classes is called the \emph{Brauer group} $\Br(X)$. Addition is given by
tensor products, and inverses come from taking dual bundles.

The Brauer group of a complex-analytic space is hard to handle. More
approachable is the \emph{cohomological Brauer group}
$\Br'(X)$, which is defined as the
torsion part of $H^2(X,\O_X^\times)$. The group extension 
$$
1\lra\O_X^\times\lra\GL_r(\O_X)\lra\PGL_r(\O_X)\lra 1 
$$ 
yields a
nonabelian coboundary map $H^1(X,\PGL_r(\CC))\ra H^2(X,\O_X)$, 
which induces an
inclusion $\Br(X)\subset \Br'(X)$.

It is possible to compute $\Br'(X)$ as an abstract group using the
exponential sequence
$0\ra\ZZ\stackrel{2\pi i}{\ra}\O_X\stackrel{\exp}{\ra}\O_X^\times\ra 1$.
The corresponding long exact sequence reads
$$
\Pic(X)\ra H^2(X,\ZZ)\ra H^2(X,\O_X)\ra H^2(X,\O_X^\times)\ra H^3(X,\ZZ)\ra H^3(X,\O_X).
$$
Let $T=T(X)$ be the torsion part of $H^3(X,\ZZ)$, which is also
the image of the coboundary map $\Br'(X)\ra H^3(X,\ZZ)$.
Furthermore, let $A=A(X)$ be the quotient $H^2(X,\ZZ)/\Pic(X)$,
which is the image of $H^2(X,\ZZ)\ra H^2(X,\O_X)$ as well.
Note that the torsion free group $A$ is sometimes called
the \emph{transcendental lattice}.

\begin{proposition}
\mylabel{computation brauer}
For any complex-analytic  space $X$, the
cohomological Brauer group canonically sits in a short exact sequence
$0\ra A\otimes\QQ/\ZZ\ra\Br'(X)\ra T\ra 0$.
\end{proposition}

\proof
We have an exact sequence
\begin{equation}
\label{long exact}
0\lra A\lra  H^2(X,\O_X)\lra H^2(X,\O_X^\times)\lra H^3(X,\ZZ)\lra H^3(X,\O_X).
\end{equation}
Set $K=H^2(X,\O_X)/A$.
Applying the functor $\Tor_i(\cdot,\QQ/\ZZ)$,
 we obtain an exact sequence
$$
0\lra \Tor_1(K,\QQ/\ZZ)\lra A\otimes\QQ/\ZZ\lra H^2(X,\O_X)\otimes\QQ/\ZZ.
$$
The term on the right vanishes, being the tensor product of
a divisible group with a torsion group.
Let $M$ be  the image of $H^2(X,\O_X^\times)\ra H^3(X,\ZZ)$.
We clearly have $T\subset M$, and $T$ equals the torsion part of $M$.
As above, we have an exact sequence
$$
0\lra\Tor_1(K,\QQ/\ZZ)\lra \Tor_1(H^2(X,\O_X^\times),\QQ/\ZZ)\lra
\Tor_1(M,\QQ/\ZZ)\lra K\otimes\QQ/\ZZ
$$
The term on the right vanishes, since $K$ is divisible and $\QQ/\ZZ$ is torsion.
The assertion now follows from the fact that $\Tor_1(M,\QQ/\ZZ)$
is the torsion part of
any abelian group $M$.
\qed

\medskip
The short exact sequence $0\ra A\otimes\QQ/\ZZ\ra\Br'(X)\ra T\ra 0$ splits,
being an extension by a divisible and hence injective group.
We call $A\otimes \QQ/\ZZ$ the \emph{analytic part} of
the cohomological Brauer group, and $T$ the \emph{topological part}.
This is a minor abuse of notation, because the short exact sequence
has no canonical splitting, but it should not cause any confusion.

For compact complex-analytic  spaces we  therefore have a  noncanonical
isomorphism
$\Br'(X)\simeq(\QQ/\ZZ)^{b_2-\rho}\oplus T$,
where $b_2$ is the second Betti number
and $\rho$ is the Picard number. The latter is defined as the rank
of the image of the coboundary map $c_1:\Pic(X)\ra H^2(X,\ZZ)$.
Hodge theory gives additional information:

\begin{corollary}
\mylabel{analytic part}
Let $X$ be a smooth compact complex-analytic space.
Assume that  $X$ is either  K\"ahler or 2-dimensional.
Then the  analytic part $A\otimes\QQ/\ZZ$ of the cohomological Brauer group $\Br'(X)$
vanishes if and only if $H^2(X,\O_X)=0$.
\end{corollary}

\proof
The condition is sufficient according to Proposition \ref{computation brauer}
and the exact sequence (\ref{long exact}).
For the converse we use the fact
from Hodge theory that the complexification of the canonical
map $H^2(X,\ZZ)\ra H^2(X,\O_X)$ is surjective.
Compare \cite{Barth; Peters; Van de Ven 1984}, Chapter IV, Proposition 2.11
for surfaces,
and \cite{Voisin 2002}, page 161 for  K\"ahler manifolds.
\qed

\medskip
Let me point out that the topological part $T$ of the Brauer group,
which is the torsion part of the cohomology group $H^3(X,\ZZ)$, is also isomorphic 
to the torsion part of the homology group $H_2(X,\ZZ)$, by the Universal Coefficient Theorem.

We next discuss the relation between algebraic and analytic theories.
The Brauer group and the cohomological Brauer group for schemes $Y$
are defined as above, but one has to use the \'etale topology instead
of the Zariski topology.
If $Y$ is an algebraic $\CC$-scheme, we have an associated complex-analytic
space $X=Y^\an$. For compact spaces, this does not influence cohomological
Brauer groups:

\begin{proposition}
\mylabel{etale classical}
Let $Y$ be an algebraic $\CC$-scheme,
and $X=Y^\an$ the associated complex-analytic space.
Then the canonical map $\Br'(Y)\ra\Br'(X)$ is surjective.
It is even bijective provided $X$ is compact.
\end{proposition}

\proof
Fix an integer $n\geq 0$. The Kummer sequence
$0\ra \mu_n\ra \O_X^\times\stackrel{n}{\ra}\O_X^\times\ra 1$
gives a short  exact sequence
$$
0\lra\Pic(X)_n\lra H^2(X,\mu_n)\lra {}_n\!\Br'(X)\lra 0.
$$
Here $\mu_n=\mu_n(\CC)$ is the group of $n$-th roots of unity,
and ${}_n\!\Br'(X)$ and $\Pic(X)_n$ are the kernel and cokernel for
multiplication-by-$n$ map.
There is a similar short exact sequence for the
\'etale topology on the
algebraic $\CC$-scheme $Y$, and we obtain a commutative diagram
$$
\begin{CD}
0 @>>> \Pic(Y)_n@>>> H^2(Y,\mu_n)@>>> {}_n\!\Br'(Y)@>>> 0\\
@. @VVV @VVV @VVV\\
0 @>>> \Pic(X)_n@>>> H^2(X,\mu_n)@>>> {}_n\!\Br'(X)@>>> 0.
\end{CD}
$$
The vertical map in the middle is bijective by
comparison results in \'etale cohomology
(\cite{SGA 4c}, Expos\'e XVI, Theorem 4.1).
It follows that $\Br'(S)\ra\Br'(X)$ is surjective.

Now suppose in addition that $X$ is compact.
It then easily follows from Serre's Gaga Theorems
(\cite{Serre 1955}, Proposition 18) that $\Pic(Y)\ra\Pic(X)$
is bijective, hence the assertion.
\qed

\medskip
Brauer groups behave similarly:

\begin{proposition}
\mylabel{etale classic'}
Notation as above. Suppose that $X=Y^\an$ is compact. Then
the canonical map $\Br(Y)\ra\Br(X)$ is bijective.
\end{proposition}

\proof
Injectivity follows directly from Proposition \ref{etale classical}.
To check surjectivity, suppose we have a holomorphic
$\PGL_r(\CC)$-principal bundle $P\ra X$. It comes
from a cocycle $\lambda_{ij}$ for some open covering $U_i\subset X$.
The problem here is that the open subsets $U_i$ are not necessarily
Zariski open, and the holomorphic maps $\lambda_{ij}\ra\PGL_r(\CC)$
are not necessarily algebraic. We sidestep the problems as follows:

The cocycle $\lambda_{ij}$ also defines, via the conjugacy action,  a holomorphic
\emph{Azumaya algebra} $\shA$ over $\O_X$. This means that $\shA$ is
a twisted form of the matrix algebra $\Mat_r(\O_X)$.
According to Serre's Gaga Theorem, the underlying
locally free $\O_X$-module is algebraic. Moreover, the structure
map $\shA\otimes\shA\ra\shA$ defining the algebra structure is algebraic
(\cite{Serre 1955}, Theorem 2).
Summing up, the Azumaya $\O_X$-algebra $\shA$ is isomorphic to
$\shB\otimes_{\O_{Y}}\O_X$ for
some locally free $\O_Y$-algebra $\shB$. Since the fibers
$\shB(x)=\shA(x)$, $x\in X$
are matrix algebras over $\CC$, the $\O_Y$-algebra $\shB$ is Azumaya.
It follows that the map on Brauer groups is surjective.
\qed

%===========================================================
\section{Nonpurity results}
\mylabel{nonpurity results}

So far we saw that analytic and algebraic Brauer groups
coincide on compact algebraic spaces.
In this section I discuss a striking difference between
analytic and algebraic Brauer groups for noncompact spaces.
For simplicity I confine the discussion to dimension two.
Throughout the paper, a \emph{surface} is a
complex-analytic space $S$ that  is
irreducible and of complex dimension two.

\begin{proposition}
\mylabel{siu vanishing}
Let $S$ be a noncompact surface.
Then $H^2(S,\O_S^\times)=H^3(S,\ZZ)$.
\end{proposition}

\proof
We have $H^2(S,\O_S)=0$ according to Siu's vanishing result for noncompact spaces
\cite{Siu 1969}, and $H^3(S,\O_S)=0$ by dimension reasons.
The exponential sequence gives an exact sequence
$$
H^2(S,\O_S)\lra H^2(S,\O_S^\times)\lra H^3(S,\ZZ)\lra H^3(S,\O_S),
$$
and the assertion follows.
\qed

\medskip
This means that the analytic part of the cohomological Brauer group
vanishes upon restrictions. The topological part behaves differently:

\begin{proposition}
Let $S$ be a smooth surface, $A\subset S$ be a discrete
subset, and $U=S-A$ the open complement.
Then the restriction map  $H^3(S,\ZZ)\ra H^3(U,\ZZ)$
induces a bijection on torsion parts.
\end{proposition}

\proof
The long exact sequence for local cohomology groups gives an exact sequence
$$
H^3_A(S,\ZZ)\lra H^3(S,\ZZ)\lra H^3(U,\ZZ)\lra H^4_A(S,\ZZ).
$$
The term on the left vanishes and the term on the right is free.
This follows from the  Thom isomorphism $H^{p-4}(A,\ZZ)\ra H^p_A(S,\ZZ)$,
compare \cite{Iversen 1986},  Chapter VIII, Proposition 2.3.
\qed

\medskip
It is more difficult to understand the behavior of
the topological part if one removes more that just points.
Recall that a complex-analytic space $X$ is called \emph{Stein}
if $H^p(X,\shF)=0$ for all $p\geq 1$ and all coherent $\O_X$-modules $\shF$.
This is the analogue of affine schemes in complex-analytic geometry.
Indeed, any complex-analytic space is covered by Stein open subsets,
and affine algebraic spaces are Stein.

\begin{proposition}
\mylabel{stein vanishing}
Suppose $S$ is a complex-analytic surface that is Stein. Then  we have $H^2(S,\O_S^\times)=0$.
\end{proposition}

\proof
Again we use the exact sequence
$$
H^2(S,\O_S)\lra H^2(S,\O_S^\times)\lra H^3(S,\ZZ).
$$
The term on the left vanishes by the Stein condition.
The term on the right also vanishes:
According to Hamm's result, any Stein space $X$ of dimension $n$
has the homotopy type of a CW-complex with cells of dimension $\leq n$ only
\cite{Hamm 1983}.
\qed

\medskip
A complex-analytic space $X$ is called \emph{holomorphically convex} if
there is a  Stein space $Y$, together with
a proper holomorphic map
$X\ra Y$. This notion is somewhere between compact spaces and Stein spaces.
It is a rather important class of spaces:
According to the \emph{Shafarevich Conjecture}, the universal covering
of any smooth projective  space should be holomorphically convex.

Note  that we may replace $Y$ by the analytic spectrum of $f_*(\O_X)$.
Then the map $f:X\ra Y$ is surjective with connected fibers, and $Y$ 
is called the \emph{Stein reduction} of $X$.

\begin{proposition}
\mylabel{holomorphically convex}
Let $S$ be a complex-analytic surface that is holo\-morphically convex
with 2-dimensional Stein reduction.
Then $H^2(S,\O_S^\times)=0$.
\end{proposition}

\proof
Let $f:S\ra Y$ be the Stein reduction.
The spectral sequence
$$
H^p(Y,R^qf_*(\O_S))\Longrightarrow H^{p+q}(S,\O_S)
$$
together with $R^2f_*(\O_S)=0$ and Steinness of $Y$ tells us that $H^2(S,\O_S)$ vanishes.
In light of the exact sequence
$$
H^2(S,\O_S)\lra H^2(S,\O_S^\times)\lra H^3(S,\ZZ),
$$
it suffices  to check that $H^3(S,\ZZ)=0$.
Consider the spectral sequence
$$
H^p(Y,R^qf_*(\ZZ))\Longrightarrow H^{p+q}(S,\ZZ).
$$
Note that the fibers $S_y=f^{-1}(y)$ are of
complex dimension $\leq 1$, and the base change maps
$R^nf_*(\ZZ)_y\ra H^n(S_y,\ZZ)$ are bijective (see \cite{Iversen 1986},
Section III, Theorem 6.2). This implies $R^qf_*(\ZZ)=0$ for $q\geq 3$.
Next we use the fact that $f:S\ra Y$ is bijective over the complement
of a discrete set $D\subset Y$. For $q\geq 1$, the sheaves $R^qf_*(\ZZ)$ are supported on
$D$, and hence $H^p(Y,R^qf_*(\ZZ))=0$ for $p,q\geq 1$.
We finally examine the terms $H^p(Y,f_*(\ZZ))$.
Note that $f_*(\ZZ)=\ZZ$, because $f:S\ra Y$ has connected fibers. According to Hamm's
result, $Y$ has the homotopy type of a CW-complex with
cells of dimension $\leq 2$ only \cite{Hamm 1983}.
This implies $H^p(Y,\ZZ)=0$ for $p\geq 3$.

Summing up, the terms $H^p(Y,R^qf_*(\ZZ))$ in the spectral sequence
vanish whenever  $p+q=3$, and therefore $H^3(S,\ZZ)=0$.
\qed

\begin{remark}
As Fabrizio Catanese pointed out to me, the
preceding result may be useful in connection with the
Shafarevich Conjecture. Suppose we want to refute the Shafarevich Conjecture.
Then we might try to find a smooth projective surface $S$ whose universal
covering $\tilde{S}$ has 2-dimensional Stein reduction,
 together with an Azumaya $\O_S$-algebra $\shA$ that does not become
a matrix algebra on  $\tilde{S}$.
I do not know whether this is feasible. But it reminds me
about the Brauer--Manin obstruction, which was used to
construct counterexamples to the Hasse principle.
\end{remark}

\medskip
Back to the comparison of algebraic and analytic Brauer groups.
Suppose that $Y$ is a smooth proper 2-dimensional $\CC$-scheme,
and $V\subset Y$ is an open subscheme. According to Grothendieck's
Purity Theorems \cite{GB II}, Section 6, the restriction map
$\Br'(Y)\ra\Br'(V)$ is injective. In contrast, we just saw
that the restriction maps on the corresponding analytic cohomological Brauer groups
is usually not injective. In other words,
algebraic and analytic Brauer groups might differ dramatically
on noncompact surfaces.

We close this section with an amusing  application to holomorphic vector bundles:
Recall that an easy computation with \v{C}ech cocycle
reveals that the restriction map of analytic Picard groups
$$
\Pic(\CC^2)\lra\Pic(\CC^2-\left\{0\right\})
$$
has infinite cokernel.
Such a behavior is rather typical:

\begin{proposition}
Let $S$ be a smooth compact algebraic surface with $H^2(S,\O_S)\neq 0$.
Fix a point $s\in S$.
Then there are infinitely many locally free coherent sheaves $\shE_U$
on the analytic surface
$U=S-\left\{s\right\}$ that do not extend to  coherent sheaves on $S$.
\end{proposition}

\proof
First note that $\Br(S)=\Br'(S)$ by Grothendieck's result
on algebraic surfaces \cite{GB II} and Proposition \ref{etale classic'}.
The analytic part $A\otimes\QQ/\ZZ$ of $\Br(S)$ is nontrivial and hence
infinite according to Corollary \ref{analytic part}.
Pick a nonzero class $\beta\in \Br(S)$ from the analytic part
and represent it by some Azumaya $\O_S$-algebra $\shA$.
The restriction $\beta_U\in\Br'(U)$ vanishes by
Proposition \ref{siu vanishing}, and this implies that  there is a
locally free $\O_U$-module $\shE_U$ with $\shA_U=\shEnd(\shE_U)$.

Suppose $\shE_U$ extends to a coherent $\O_S$-module $\shE$.
Passing to double duals, we may assume that $\shE$ is locally free.
The bijection $\shEnd(\shE)_U\ra\shA_U$ extends
to a map $\shEnd(\shE)\ra \shA$. Its determinant vanishes either nowhere or in
 codimension one, and we infer that the map is bijective.
A similar argument shows that this map is an isomorphism of algebras.
This implies $\beta=0$, contradiction.
\qed

%===========================================================
\section{Serre's good groups}
\mylabel{serre's good groups}

Serre's notion of good groups plays a crucial role in the sequel, and
I want to recall this concept now.
Let $G$ be a group and $\hat{G}=\invlim G/N$ be its profinite
completion. Here the inverse limit runs over all normal subgroups
$N\subset G$
of finite index. We regard both $G$ and $\hat{G}$ as topological groups:
the group $G$ carries the discrete topology, and $\hat{G}$
is endowed with the  inverse limit topology where the factors $G/N$ are
discrete.

Now let $M$ be a finite $G$-module. By this I understand a finite abelian
group with discrete topology, and having a $G$-module structure.
The map $G\ra\Aut(M)$ factors over $\hat{G}$,  and we may regard $M$
as a topological $\hat{G}$-module as well.
The canonical map $G\ra\hat{G}$ induces a restriction map
$$
H^p(G,M)\lra H^p(\hat{G},M),\quad p\geq 0
$$
on cohomology groups.
Note that our cohomology groups are defined in terms of continuous
cochains.
Serre showed in \cite{Serre 1972}, Chapter I, \S 2.6
that these restriction maps
are bijective for $p=0,1$.
He calls a group $G$ \emph{good} if the restriction maps
$H^p(G,M)\ra H^p(\hat{G},M)$ are bijective for
all integers $p\geq 0$ and all finite $G$-modules $M$.
Finite groups are clearly good groups.
Recall that a group is called \emph{almost free} if
it contains a free subgroup of finite index.

\begin{proposition}
\mylabel{free group}
Almost free groups are good groups.
\end{proposition}

\proof
First consider the case that $G$ is a  free group. Then $H^p(G,M)=0$ for $p\geq 2$ and
any $G$-module $M$, because $G$ is the fundamental group of
an Eilenberg--Maclane space $K(G,1)$ with cells of dimension $\leq 1$ only.
On the other hand, we have $H^p(\hat{G}, M)=0$ for $p\geq  2$ and any
torsion $G$-module $M$ by \cite{Serre 1972}, Chapter I, Proposition 16.

Now suppose that $G$ is almost free. Then we find a normal subgroup
$N\subset G$ that is free and with $Q=G/N$ finite.
According to \cite{Schneebeli 1979}, Lemma in Section 5.1,
the canonical map on profinite completions $\hat{N}\ra\hat{G}$ is
injective with $Q=\hat{G}/\hat{N}$.

The following is a variant of an argument due to Serre (\cite{Serre 1972}, Section 2.6):
Let $M$ be a finite $G$-module, and consider the
two Hochschild--Serre spectral sequences
$$
H^p(Q,H^q(N,M))\Longrightarrow H^{p+q}(G,M),\quad\quad
H^p(Q,H^q(\hat{N},M))\Longrightarrow H^{p+q}(\hat{G},M).
$$
Using that $N$ is free we obtain a long exact sequences
$$
\ldots \ra  H^p(Q,H^0(N)) \ra  H^p(G) \ra H^{p-1}(Q,H^1(N)) \ra H^{p+1}(Q,H^0(N))
\ra\ldots
$$
and another long exact sequence
$$
\ldots \ra  H^p(Q,H^0(\hat{N})) \ra  H^p(\hat{G}) \ra H^{p-1}(Q,H^1(\hat{N})) \ra H^{p+1}(Q,H^0(\hat{N}))
\ra\ldots,
$$
where the coefficient groups are always $M$.
Note that the canonical mappings
$H^i(N,M)\ra H^i(\hat{N},M)$ are always bijective for $i=0,1$.
Using the 5-lemma, we deduce that the canonical map $H^p(G,M)\ra H^p(\hat{G},M)$
is bijective.
\qed

\medskip
Recall that a group $G$ is called \emph{polycyclic}
if there is a finite sequence of subgroups
$0=G_0\subset G_1\subset\ldots\subset G_n=G$ so that
$G_{i-1}\subset G_i$ are normal with cyclic factors $G_i/G_{i-1}$.
This are precisely the solvable groups all whose subgroups are finitely generated.
Note that  finitely generated nilpotent groups are polycyclic.
A group is called \emph{almost polycyclic} if it contains
a polycyclic subgroup of finite index.

\begin{proposition}
\mylabel{good group}
Almost polycyclic  groups are good groups.
\end{proposition}

\proof
Suppose first that $G$ is a polycyclic group.
Then the cohomology groups $H^p(G,M)$ are finite for
all finite $G$-modules $M$ and all integers $p\geq 0$.
This is obvious for $G$ cyclic, and follows by induction
on the length $n$ of the subnormal series $G_i\subset G$,
together with the Hochschild--Serre spectral sequence.

To proceed, we use a general result of Serre:
Let $G$ be an arbitrary group containing a normal subgroup $N\subset G$.
Suppose that $N$ and $Q=G/N$ are good, and that $H^p(N,M)$ are finite
for all finite $A$-modules $M$ and $p\geq 0$.
Then Serre outlined in \cite{Serre 1972}, Section 2.6 that this implies
that $G$ is also good: He first checks that the sequence
$1\ra\hat{N}\ra\hat{G}\ra\hat{Q}\ra 1$ remains a group extension,
and then compares the two Hochschild--Serre spectral sequences.
Using this, we easily verify that our polycyclic group $G$ is good, again by
induction on the length $n$ of the subnormal series $G_i\subset G$.

Now suppose that $G$ is almost polycyclic. Then we find a
normal subgroup $N\subset G$ that is polycyclic and with $Q=G/N$ finite.
Then $N,Q$ are good, and $H^q(N,M)$ are finite for all finite $G$-modules $M$.
Repeating Serre's argument as above, we see that $G$ is good.
\qed

\medskip
Good groups are very useful with respect to Grothendieck's question
on Brauer groups.
Let $G$ be a group. The exact sequence of groups with trivial $G$-action
$$
1\lra\CC^\times\lra\GL_r(\CC)\lra\PGL_r(\CC)\lra 1
$$
induces a coboundary map in nonabelian cohomology
\begin{equation}
\label{coboundary map}
H^1(G,\PGL_r(\CC))\lra H^2(G,\CC^\times)
\end{equation}
as explained in \cite{Serre 1972}, Chapter I, \S 5.
In representation theory of finite groups, $H^2(G,\CC^\times)$ is also
called the \emph{Schur multiplier} (confer \cite{Karpilovsky 1987}).
In this context, Grothen\-dieck's question is: Given a torsion class
$\beta\in H^2(G,\CC^\times)$, does there exist some $r>0$
so that $\beta$ lies in the image of the coboundary map
(\ref{coboundary map})?

\begin{proposition}
\mylabel{coboundary good}
Let $G$ be a good group. Then any torsion class $\beta\in H^2(G,\CC^\times)$
lies in the image of the coboundary map
$H^1(G,\PGL_r(\CC))\ra H^2(G,\CC^\times)$ for some integer $r>0$.
\end{proposition}

\proof
The idea is to reduce to the case that the group $G$ is finite.
Let $\mu_n=\mu_n(\CC)$ be the group of complex $n$-th roots of unity.
Using the Kummer sequence, we infer that $\beta$ lies in the
image of $H^2(G,\mu_n)\ra H^2(G,\CC^\times)$, whenever $n\beta=0$.
Choose $\alpha\in H^2(G,\mu_n)$ mapping to $\beta$.

Since $G$ is good, we have
$$
H^2(G,\mu_n) = H^2(\hat{G},\mu_n) =\dirlim H^2(G/N,\mu_n).
$$
Replacing $G$ by some suitable $G/N$, we may assume that
the group $G$ is finite.
It is then a fact from representation theory that any
factor system comes from a projective representation
(see for example \cite{Huppert 1967}, Chapter V, Hilfssatz 24.2),
hence the result.
\qed

%===========================================================
\section{Applications to complex spaces}
\mylabel{applications to complex spaces}

We now come back to geometry.
Fix a complex-analytic space $X$.
We seek conditions under which the canonical inclusion
$\Br(X)\subset\Br'(X)$ is an equality.
The following result is interesting because it makes only assumptions
on the homotopy type of $X$:

\begin{theorem}
\mylabel{brauer good}
Let $X$ be a  connected complex-analytic space.
Suppose that the fundamental group $\pi_1(X)$ is good, and that the
subgroup of $\pi_1(X)$-invariants in  the Pontrjagin dual $\Hom(\pi_2(X),\QQ/\ZZ)$ vanishes.
Then $\Br(X)=\Br'(X)$.
\end{theorem}

\proof
Let $\beta\in\Br'(X)$, say with $n\beta=0$, and choose a class
$\alpha\in H^2(X,\mu_n)$ mapping to $\beta$.
Let $\tilde{X}\ra X$ be the universal covering.
Then $\pi_2(X)=\pi_2(\tilde{X})$, and the canonical map
$H_2(\tilde{X},\ZZ)\ra \pi_2(\tilde{X})$ is bijective by
the Hurewicz Theorem. Moreover, the canonical
map $H^2(\tilde{X},\mu_n)\ra \Hom(H_2(\tilde{X},\ZZ),\mu_n)$
is bijective by the Universal Coefficient Theorem.
The universal covering yields a spectral sequence
$$
H^p(G,H^q(\tilde{X},\mu_n))\Longrightarrow H^{p+q}(X,\mu_n),
$$
where $G=\pi_1(X)$. Since $H^1(\tilde{X},\mu_n)=0$,
this yields  a short exact sequence
$$
0\lra H^2(G,\mu_n)\lra H^2(X,\mu_n)\lra H^0(G,H^2(\tilde{X},\mu_n)).
$$
The term on the right vanishes, in light of our assumption on $\pi_2(X)$
and the  inclusion
$H^2(\tilde{X},\mu_n)\subset\Hom(\pi_2(X),\QQ/\ZZ)$.
The upshot is that any 2-cohomology class on $X$ with values in $\mu_n$
comes from group cohomology.

To finish the proof we use the assumption that $G=\pi_1(X)$ is good, which
implies
$$
H^2(G,\mu_n)=H^2(\hat{G},\mu_n)=\dirlim H^2(G/N,\mu_n).
$$
Hence there is an finite \'etale covering
$g:X'\ra X$ with $g^*(\beta)=0$, and in particular
$g^*(\alpha)=0$. Then  a result
of Hoobler  (\cite{Hoobler 1972}, Section 3), see also Gabber (\cite{Gabber 1981}, Lemma 4), tells us that
$\alpha\in\Br(X)$.
Alternatively, we could construct a principal $\PGL_r(\CC)$-bundle
$P\ra X$ representing $\beta$
as the principal bundle associated with
the representation $\pi_1(X)\ra\PGL_r(\CC)$
constructed in the end of the proof for Proposition \ref{coboundary good}.
\qed

\medskip
The preceding Theorem applies in particular if $\pi_2(X)=H_2(\tilde{X},\ZZ)=0$.
This is made to measure for complex Lie groups:

\begin{corollary}
\mylabel{lie brauer}
Let $X$ be any complex Lie group.
Then $\Br(X)=\Br'(X)$.
\end{corollary}

\proof
We may assume that $X$ is connected.
The fundamental group $\pi_1(X)$ is a finitely generated
abelian group, hence good.
Indeed: it is abelian because $X$ is a $H$-space, and
it is finitely generated because $X$ is homotopy equivalent
to a compact differentiable manifold (see \cite{Iwasawa 1949}, Theorem 6).
The homotopy group $\pi_2(X)$ vanishes, according to
Browder results on torsion in homology of $H$-spaces (\cite{Browder 1961},
Theorem 6.11).
Consequently Theorem \ref{brauer good} applies.
\qed

\medskip
This generalizes a result of Iversen on characterfree
linear algebraic groups \cite{Iversen 1976}.
It also generalizes results of
Hoobler \cite{Hoobler 1972} and Berkovi\v c \cite{Berkovic 1972}
on abelian varieties, and of Elencwajg and Narasimhan
on complex tori \cite{Elencwajg; Narasimhan 1983}.

It is instructive to look at the case of connected abelian Lie groups.
They have the form $X=\CC^n/\Gamma$ for some lattice
$\Gamma\subset\CC^n$, and are studied in the book of Abe and Kopfermann
\cite{Abe; Kopfermann 2001}. We have $H^i(X,\ZZ)=\Hom(\Lambda^i\Lambda,\ZZ)$,
and from this it is in principle possible to compute
the Brauer group $\Br(X)=\Br'(X)$
as a quotient of $H^2(X,\QQ/\ZZ)$.
The full cohomology group $H^2(X,\O_X^\times)$, however, can be
tremendously large: As explained in \cite{Kazama 1984},
there are lattices $\Gamma$ such that $H^2(X,\O_X)$ are infinite
dimensional complex vector spaces.

We next apply our result to Hopf manifolds.
Recall that a complex space $X$ is called a
\emph{Hopf manifold} if it is compact and its universal
covering $\tilde{X}$ is biholomorphic to $\CC^n-\left\{0\right\}$
with $n\geq 2$.
Hopf manifolds are the  simplest examples of compact complex manifolds
that are not K\"ahler.

\begin{corollary}
\mylabel{hopf brauer}
Let $X$ be a Hopf manifold. Then $\Br(X)=\Br'(X)$.
\end{corollary}

\proof
We obviously  have $H_2(\tilde{X},\ZZ)=0$.
Furthermore, Kodaira proved in \cite{Kodaira 1966}, Theorem 30
that the fundamental group sits in
a central  extension
\begin{equation}
\label{hopf fundamental}
0\lra \ZZ\lra \pi_1(X)\lra G\lra 1
\end{equation}
for some finite group $G$. Actually,
Kodaira treated the 2-dimensional case,
but his arguments work in all dimensions $n\geq 2$.
So $\pi_1(X)$ is almost free hence good by  Proposition \ref{free group},
and we conclude with Theorem \ref{brauer good}.
\qed

\medskip
Actually, we may compute the Brauer group of a $n$-dimensional
Hopf manifold $X$ as above:
Let $c\in H^2(G,\ZZ)$ be the extension class of the
central extension (\ref{hopf fundamental}).
It defines   homomorphisms $c:H^i(G,\ZZ)\ra H^{i+2}(G,\ZZ)$
via the cup product.

\begin{proposition}
\mylabel{computation hopf}
Let $X$ be a Hopf manifold. With the preceding notation,
we have an exact sequence
$H^1(G,\ZZ)\stackrel{c}{\ra} H^3(G,\ZZ)\ra \Br(X)\ra H^2(G,\ZZ)\stackrel{c}{\ra} H^4(G,\ZZ)$.
\end{proposition}

\proof
First note that the analytic part of the Brauer group vanishes,
because $H^2(X,\O_X)=0$. The latter  is due to Kodaira for Hopf surfaces
(\cite{Kodaira 1966}, Theorem 26), and follows for higher dimensional Hopf manifolds
from a result of Mall (\cite{Mall 1991}, Theorem 3).

We have to compute the topological part of the Brauer group.
To simplify notation, set $\pi_1=\pi_1(X)$ and $m=2n-1$. Note that $n\geq 2$ and $m\geq 3$.
The spectral sequence
$H^p(\pi_1,H^q(\tilde{X},\ZZ))\Rightarrow H^{p+q}(X,\ZZ)$ reduces to a
long exact sequence
$$
\ra H^{p-1-m}(\pi_1,H^{m}(\tilde{X}))\ra H^p(\pi_1,H^0(\tilde{X}))\ra
H^p(X)\ra H^{p-m}(\pi_1,H^{m}(\tilde{X}))\ra.
$$
In particular, the edge maps of
the spectral sequence
$H^p(\pi_1,H^0(\tilde{X},\ZZ))\ra H^p(X,\ZZ)$ are bijective
for $p<m$. For $p=m$, we obtain an exact sequence
$$
0\lra H^m(\pi_1,H^0(\tilde{X}))\lra H^m(X)\lra H^0(\pi_1,H^m(\tilde{X})).
$$
The term on the right equals $H^m(\tilde{X})^{\pi_1} = \ZZ$,
because $H^m(\tilde{X})=\ZZ$ and the action of $\pi_1$ on $\tilde{X}$ is
orientation preserving. Hence the torsion in
$H^m(X)$ is contained in $H^m(\pi_1,H^0(\tilde{X}))$.

To proceed, we view $\ZZ=H^m(\tilde{X})$ as a $\pi_1$-module with trivial action.
The 
Hochschild--Serre spectral sequence $H^p(G,H^q(\ZZ,\ZZ))\Rightarrow H^{p+q}(\pi_1,\ZZ)$
for the central extension
(\ref{hopf fundamental}) reduces to a long exact sequence
$$
H^{p-2}(G,H^{1}(\ZZ,\ZZ))\ra H^p(G,H^0(\ZZ,\ZZ))\ra H^p(\pi_1,\ZZ)\ra
H^{p-1}(G,H^{1}(\ZZ,\ZZ)),
$$
because $H^p(\ZZ, M)=0$ for $p\geq 2$ and any $\ZZ$-module $M$.
Since $H^p(G,M)$ are torsion groups for $p>0$ and any $G$-module $M$,
the groups $H^p(\pi_1,\ZZ)$ are torsion for 
$p\geq 2$.

Combining the preceding two paragraphs, we infer that
$H^3(\pi_1,\ZZ)$ is torsion, and equals the torsion part in $H^3(X,\ZZ)$.
Using Proposition \ref{computation brauer}, we infer
$$
\Br(X)=\Br'(X)=H^3(\pi_1,\ZZ).
$$
By the Hochschild--Serre spectral sequence, this group sits in an exact sequence
\begin{equation}
\label{cohomology hopf}
H^1(G,\ZZ)\lra H^3(G,\ZZ)\lra H^3(\pi_1,\ZZ) \lra H^2(G,\ZZ)\lra H^4(G,\ZZ).
\end{equation}
According to \cite{Hochschild; Serre 1953}, Section 6 the outer maps are
taking cup products with the extension class $-c\in H^2(G,\ZZ)$,
hence the assertion.
\qed

\medskip
To my knowledge, there has been no attempt to classify such group extensions
occurring in Hopf manifolds, except for dimension two:
Kato determined
all possible groups $\pi_1(S)$ for Hopf surfaces $S$ in
\cite{Kato 1975}, \cite{Kato 1989}.
There are two cases: In the \emph{linear case}
$\pi_1(S)$ is conjugate to a subgroup
of $\GL_2(\CC)$, and then
$G$ is an extension
of a subgroup $K\subset\SL_2(\CC)$ by a cyclic group.
Such subgroups are classified (cyclic, dihedral, tetrahedral, octahedral,
icosahedral).
In the nonlinear case $\pi_1(S)$ is not conjugate to a subgroup
of $\GL_2(\CC)$, and $G$ must be cyclic.

%===========================================================
\section{From open to compact surfaces}
\mylabel{from open to compact surfaces}

In order to tackle smooth compact surfaces, we have
to improve Theorem \ref{brauer good}.
The key result is the following, which works in dimension 2 only:

\begin{theorem}
\mylabel{zariski good} Let $S$ be a smooth compact
surface, and $V\subset S$ be a nonempty Zariski open subset.
Suppose that the fundamental group $\pi_1(V)$ is good, and that the
subgroup of $\pi_1(V)$-invariants in the Pontrjagin dual
$\Hom(\pi_2(V),\QQ/\ZZ)$ vanishes. Then $\Br(S)=\Br'(S)$.
\end{theorem}

\proof
Fix a class $\beta\in\Br'(S)$, say with $n\beta=0$, and choose
a class $\alpha\in H^2(S,\mu_n)$ mapping to $\beta$.
As in the proof of Theorem \ref{brauer good},
there is a finite \'etale covering $U\ra V$ with $\alpha_{U}=0$.
By \cite{Detloff; Grauert 1994}, Theorem 3.4
there is a compactification $U\subset X$ with some compact complex-analytic
surface $X$,
so that we have a commutative diagram
$$
\begin{CD}
U @>>> X\\
@VVV @VVV\\
V@>>> S.
\end{CD}
$$
Replacing $X$ by a suitable blowing-up with center in the 
boundary  $C=X-U$,
we may assume that the compact surface $X$ is smooth, and
that $C=C_1+\ldots +C_m$ is a normal crossing divisor with
smooth irreducible components
(see, for example, \cite{Laufer 1971}, Chapter I and II).

In the next step we use local cohomology groups $H^2_C(X,\mu_n)$. 
Combining the Kummer sequence with the local cohomology sequence on $X$,
we obtain a diagram
$$
\begin{CD}
@.@. H^2_C(X,\mu_n)\\
@.@. @VVV\\
0@>>>\Pic(X)_n @>>> H^2(X,\mu_n) @>>> {}_n\!\Br'(X)@>>>0\\
@.@. @VVV\\
@.@. H^2(U,\mu_n)
\end{CD}.
$$
Here ${}_n\!\Br'(X)$ and $\Pic(X)_n$ are the kernel and cokernel for
multiplication-by-$n$ map.
We remark that that the upper vertical map factors over the left
horizontal map.
To see this, note that according to 
\cite{Iversen 1986}, Chapter 10, Proposition 2.5, the 
canonical map between local cohomology groups
$$
\bigoplus_{i=1}^m H^2_{C_i}(X,\mu_n)\lra H^2_C(X,\mu_n)
$$
is bijective. The summands $H^2_{C_i}(X,\mu_n)$ are
freely generated by the cycle class $\cl_X(C_i)$,
as explained in \cite{Deligne 1977}, Proposition 2.2.6 on page 141
(this class is called the \emph{Thom class} in topology).
Moreover, the image of $\cl_X(C_i)$ and $\O_X(C_i)$ in
$H^2(X,\mu_n)$ coincide by \cite{Freitag; Kiehl 1988}, Chapter II, Proposition
2.2.
The upshot is that $\alpha_U=0$ implies
that $\alpha_{X}\equiv 0$ modulo $\Pic(X)_n$, and hence $\beta_{X}=0$.
To finish the proof, it therefore suffices to verify the following
statement:

\begin{lemma}
\mylabel{covering}
Let $S$ be a smooth complex-analytic surface, $X$ a normal complex-analytic surface, and
$f:X\ra S$ be a proper holomorphic surjection.
Then $\Br(S)$ contains
the kernel of the induced  map $\Br'(S)\ra\Br'(X)$.
\end{lemma}

\proof
First, consider the special case that the proper holomorphic map
$X\ra S$ is  finite.  It then follows that it is flat as well,
because $S$ is smooth and $X$ is Cohen--Macauley
(this follows from \cite{Serre 1965}, Chapter IV, Theorem 9).
Then the result of Hoobler (\cite{Hoobler 1972}, Section 3) and Gabber (\cite{Gabber 1981}, Lemma 4) tells us that
any cohomology class $\beta\in\Br'(S)$ with $\beta_X=0$ lies in 
$\Br(S)$.

Now consider the general case. There exist a proper bimeromorphic mapping
$f:S'\ra S$ so that the 2-dimensional integral  component
$Y'\subset X\times_S S'$ is flat over $S'$.
This is a special case of Hironaka's general result on flattening
\cite{Hironaka 1975}. Here we need only the 2-dimensional case, which
also  Maurer \cite{Maurer 1978} worked out.
Actually  we only need that  $Y'\ra S'$ is finite. Replacing $S'$ by a suitable
blowing-up, we may assume that $S'$ is smooth.
The normalization $X'$ of $Y'$ is finite over $S'$ as well.
According to the special case, the preimage $\beta_{S'}$ lies
in $\Br(S')$. The statement now follows from the following observation:

\begin{lemma}
\mylabel{bimeromorphic}
Let $S,S'$ be two smooth surfaces, and $f:S'\ra S$ a proper bimeromorphic
map. Then any class 
$\beta\in\Br'(S)$ with $\beta_{S'}\in\Br(S')$ lies in $\Br(S)$.
\end{lemma}

\proof
As in the proof for Proposition \ref{etale classic'}, it is more convenient
to work with Azumaya algebras than with principal bundles.
Let $\shB$ be an Azumaya $\O_{S'}$-algebra representing $\beta_{S'}$.
Then the double dual $\shA=f_*(\shB)^{\vee\vee}$ is coherent reflexive
$\O_S$-algebra, which is actually locally free, because any reflexive sheaf
on a complex manifold is locally free in codimension two. I claim that
$\shA$ is Azumaya. For this, consider the $\O_X$-algebra homomorphism
$$
\varphi:\shA\otimes\shA^\op\lra\shEnd(\shA),\quad A\otimes A'\longmapsto (B\mapsto ABA').
$$
Both sheaves are locally free, and the map is bijective on
$S-T$, where $T\subset S$ is the discrete subset of all points
$s\in S$ with 1-dimensional fiber $f^{-1}(s)$.
The set of points where the map is not bijective is a Cartier divisor
(choose local holomorphic bases and take determinants)
hence our map $\varphi$ is  bijective everywhere on $S$.
This means  that $\shA$ is an Azumaya algebra
by  \cite{GB I}, Theorem 5.

It remains to check that the Azumaya algebra $\shA$ 
represents the cohomology class $\beta$.
This is  easy for algebraic surfaces with Grothendieck's purity
results from \cite{GB III}, Section 6. We already saw, however,
that such purity results do not hold true in the analytic situation.
We shall use methods
from \v{C}ech cohomology instead. The argument runs as follows:

By Lemma \ref{blowing up} below, it suffices to check that the
Azumaya $\O_{S'}$-algebras  $\shA'=f^*(\shA)$ and
$\shB$ have the same cohomology class. We check this by construction
explicit 2-cocycles. To do this, we first have to settle some
technical points regarding the existence of 2-cocycles for 2-cohomology classes.
Choose an open covering $U_i\subset S$, $i\in I$ so that $\beta_{U_i}=0$.
After passing to a refinement, we may assume that every  $s\in T$ 
is contained in precisely one $U_i$.
The Azumaya $\O_{S'}$-algebra $\shB$ represents
the zero class on each preimage $U'_i=f^{-1}(U_i)$, hence there
is a locally free $\O_{U'_i}$-module $\shE_i$ and an isomorphism
$s_i:\shB_{U'_i}\ra\shEnd(\shE_i)$.
The double duals $\shF_i=f_*(\shE_i)^{\vee\vee}$ are locally free
$\O_{U_i}$-modules, and the bijections $s_i$ induce bijections 
$t_i:\shA_{U_i}\ra \shEnd(\shF_i)$, by the same argument as above.

On the overlaps $U_{ij}=U_i\cap U_j$ there are  invertible sheaves $\shL_{ij}$
defined by the condition $\shF_i|_{U_{ij}}\otimes\shL_{ij}\simeq\shF_{j}|_{U_{ij}}$.
We regard $U_{ij}\mapsto \shL_{ij}$ as a 1-cochain with respect to the
open covering $\foU=(U_i)_{i\in I}$ taking values in the presheaf
$\shH^1\O_S^\times$. 
The latter is defined by $\Gamma(U,\shH^1\O_S^\times)=\Pic(U)$. 
As explained in \cite{Schroeer 2003}, Lemma 3.1 the cochain $\shL_{ij}$ is actually
a cocycle.
Things would be particularly nice if the $\shL_{ij}\in\Pic(U_{ij})$ are trivial.
We can achieve this, for example, by first refining the covering to  a covering with
$U_{ij}$ contractible
(compare \cite{Bott; Tu 1986}, Corollary I.5.2), and then further refining to a covering that is Stein.
However, the following argument works in more general situations as well:
By definition, complex-analytic spaces
have a countable basis for the topology, hence are paracompact.
This implies $\cH^1(X,\shH^1\O_X^\times)=0$, as explained in \cite{Godement 1964},
Chapter II, Proposition 5.10.1.
So after passing to a refinement, we find invertible $\O_{U_i}$-modules $\shL_i$
with $\shL_i|_{U_{ij}}\otimes\shL_{ij}\simeq\shL_{j}|_{U_{ij}}$. 
Replacing $\shE_i$ by $\shE_i\otimes f^*(\shL_i)$, which replaces
$\shF_i$ by $\shF_i\otimes \shL_i$, we may assume 
that all $\shL_{ij}$ are trivial.

Summing up, we placed ourselves into a situation in which there are
isomorphisms $t_{ij}:\shF_j|_{U_{ij}}\ra\shF_i|_{U_{ij}}$.
Choose such bijections subject to the conditions $t_{ii}=\id$ and $t_{ij}t_{ji}=\id$.
Now consider the  2-cochain $\lambda\in Z^2(\foU,\O_S^\times)$ defined by
$$
t_{jk}t_{ij}=\lambda_{ijk}t_{ik}
$$
on the triple overlaps $U_{ijk}=U_i\cap U_j\cap U_k$.
According to Giraud (\cite{Giraud 1971}, Chapter IV, Section 3.5, see also
Chapter V, Section 4)
this is a 2-cocycle whose cohomology class
represents the cohomology class of $\shA$.

By our construction, the holomorphic map $U'_{ij}\ra U_{ij}$ is biholomorphic
whenever $i\neq j$. For such pair of indices, the
isomorphisms $t_{ij}$ may be regarded as isomorphisms
$t_{ij}':\shE_i|_{U'_{ij}}\ra\shE_j|_{U'_{ij}}$.
In case $i=j$ we define $t'_{ii}=\id$. This yields another 2-cocycle
$\lambda'\in Z^2(\foU',\O_{S'}^\times)$ via
$$
t'_{jk}t'_{ij}=\lambda'_{ijk}t'_{ik},
$$
whose cohomology class represents the cohomology class of $\shB$.
By construction, the local sections $\lambda_{ijk}\in\Gamma(U_{ijk},\O_S)$
map to the local section $\lambda'_{ijk}\in\Gamma(U'_{ijk},\O_{S'})$ under
the canonical map. It follows
that the Azumaya $\O_{S'}$-algebras $\shA'$ and $\shB$
have the same class.
\qed

\medskip
It remains to check the following observation:
\begin{lemma}
\mylabel{blowing up}
Let $f:S'\ra S$ be a proper bimeromorphic map between smooth surfaces.
Then the canonical map  $H^2(S,\O_S^\times)\ra H^2(S',\O_{S'}^\times)$ is bijective.
\end{lemma}

\proof
Compare the long exact sequence
$$
\Pic(S)\lra H^2(S,\ZZ)\lra H^2(S,\O_S)\lra H^2(S,\O_S^\times)\lra H^3(S,\ZZ)
$$
with the corresponding sequence for $S'$, and use
the vanishing $R^pf_*(\O_{S'})=0$ for $p>0$ and
$R^pf_*(\ZZ)=0$ for $p\neq 0,2$.
\qed

%===========================================================
\section{Elliptic surfaces}
\mylabel{elliptic surfaces}

In this section we apply Theorem \ref{zariski good} to elliptic surfaces.
Recall that an \emph{elliptic surface} is a smooth compact surface $S$,
together with a holomorphic map $f:S\ra B$ onto a smooth compact curve $B$
so that almost all fibers $S_b=f^{-1}(b)$, $b\in B$ are elliptic
curves.
Elliptic surfaces might have any Kodaira dimension
$\kappa(S)\leq 1$, and any algebraic dimension $a(S)\leq 2$.
Recall that the \emph{algebraic dimension} $a(S)$ is the
transcendence degree of the field of meromorphic functions on $S$.
Surfaces with $a(S)=2$ are algebraic, surfaces with $a(S)=1$ are elliptic,
and surfaces with $a(S)=0$ contain only finitely many curves.
An algebraic surface might have several elliptic structures,
even infinitely many. In contrast, surfaces with $a(S)\leq 1$ have
at most one elliptic structure, if any.

\begin{proposition}
\mylabel{elliptic brauer} Let $S$ be an elliptic surface. Then
$\Br(S)=\Br'(S)$.
\end{proposition}

\proof
Fix an elliptic fibration $f:S\ra B$. It is possible to
compute the fundamental group $\pi_1(S)$ in terms of this
fibration, see \cite{Xiao 1991}. This fundamental group, however,
might be rather small, due to the influence of singular fibers.
Things simplify much if one throws away the singular fibers:

Choose finitely many points $b_1,\ldots,b_m\in B$ such that all
singular fibers occur among the fibers $S_{b_1},\ldots,S_{b_m}$,
and that the complement $V=B-\left\{b_1,\ldots b_m\right\}$ is
\emph{hyperbolic}. The latter means that its universal covering
space is the upper half plane $\tilde{V}=\foH$. Set
$U=S-\bigcup_{i=1}^m S_{b_i}$. Then the induced map $U\ra V$ is
proper and smooth with elliptic fibers. Let $F\subset U$ be any
fiber. The long homotopy sequence reads $$
\pi_2(V)\lra\pi_1(F)\lra \pi_1(U)\lra \pi_1(V)\lra 0. $$ By
construction, the group $\pi_1(V)$ is free and
$\pi_2(V)=\pi_2(\tilde{V})$ vanishes. We conclude that $\pi_1(U)$
is an extension of a free  group by $\ZZ^{\oplus 2}$.
Arguing as in the proof for Proposition \ref{good group},
we see that $\pi_1(U)$ is good.
Moreover, the universal covering $\tilde{U}=\foH\times\CC$ is
contractible. This means that Proposition \ref{zariski good}
applies to our situation, and we conclude $\Br(S)=\Br'(S)$.
\qed

\medskip
Let me close this section with some remarks  about the relation of
the elliptic surface $S\ra B$ to the corresponding jacobian
elliptic surface $X\ra B$.
According to Kodaira \cite{Kodaira 1963}, any elliptic surface $f:S\ra B$
(without exceptional curves and multiple fibers) comes along with two invariants: The \emph{homological invariant} $\shG=R^1f_*(\ZZ)$, and
the \emph{functional invariant} $j:B\ra\PP^1$. The latter attaches
to each point with smooth fiber the $j$-invariant of its fiber.
Kodaira showed that the jacobian fibration $X\ra B$ of $S\ra B$ has
the same homological and functional invariant.
By construction, the jacobian fibration has a section, so $X$ is an algebraic
surface. It is not difficult to see that the topological part of $\Br(S)$ and
$\Br(X)$ coincide. In contrast, the analytic part might differ drastically due to jumps in  Picard numbers.
The situation is simpler for algebraic surfaces $S$: Here Nori \cite{Nori 1986}
showed that $\rho(S)=\rho(X)$, and hence $\Br(S)\simeq\Br(X)$.

%===========================================================
\section{Surfaces of class VII}
\mylabel{Surfaces of class VII}

In this section we analyze smooth compact surfaces $S$ of class VII.
By definition, this means $b_1(S)=1$.
Such surfaces are not algebraic, not even K\"ahler.
They form the only class of surfaces resisting complete classification.
We first observe that the cohomological Brauer group is rather
small:

\begin{lemma}
\mylabel{small}
For surfaces $S$ of class VII we have $H^2(S,\O_S^\times)= H^3(S,\ZZ)$
\end{lemma}

\proof
The exponential sequence gives an exact sequence
$$
H^2(S,\O_S)\lra H^2(S,\O_S^\times)\lra H^3(S,\ZZ)\lra H^3(S,\O_S).
$$
The term on the right vanishes for dimension reason.
and the term on the left also vanishes, as Kodaira showed in \cite{Kodaira 1966},
Theorem 26.
\qed

\medskip
Recall that a \emph{global spherical shell} consists of
an holomorphic open embedding of some
$U=\left\{z\in\CC^2\mid 1-\epsilon<|z|<1+\epsilon\right\}$ with $0<\epsilon<1$
into $S$ so that the complement $S-U$ is connected.
A surface admitting a global spherical shell is of class VII.
Such surfaces are not very interesting with respect to Grothendieck's question on
Brauer groups:

\begin{proposition}
\mylabel{shell brauer}
Let $S$ be a smooth compact surface containing
a global spherical shell. Then $\Br(S)=\Br'(S)=0$.
\end{proposition}

\proof
We just saw that the analytic part of $\Br'(S)$ vanishes.
Using a Mayer--Vietoris argument,
Dloussky showed in \cite{Dloussky 1984}, Lemma 1.10
that $H_2(S,\ZZ)$ is torsion
free. In other words, $H^3(S,\ZZ)$ is torsion free,
so the topological part of the Brauer group vanishes
as well.
\qed

\medskip
We now turn to surfaces of class VII with $b_2=0$.
Such surfaces are indeed classified.
If there is a curve $C\subset S$, then Kodaira proved that
$S$ is a Hopf surface \cite{Kodaira 1966}, Theorem 34.
We already settled this case in Corollary \ref{hopf brauer}.
If there is no curve at all on $S$, then the classification results
of Inoue apply \cite{Inoue 1974}.
He showed that the universal covering is $\tilde{S}=\foH\times\CC$
and that $\pi_1(S)$ is polycyclic.
Actually, Inoue made an additional technical assumption,
namely the existence of a twisted vector field.
Later, the work of Bogomolov \cite{Bogomolov 1976}, Li, Yau, and
Zheng \cite{Li; Yau; Zheng 1990}, \cite{Li; Yau; Zheng 1994}, and
Teleman \cite{Teleman 1994} revealed that this assumption automatically holds.
One refers to surfaces with $b_1=1$ and $b_2=0$ without curves as
\emph{Inoue surfaces}.
Summing up, we may apply Theorem \ref{brauer good} and deduce:

\begin{proposition}
\mylabel{hopf inoue}
Let $S$ be a smooth compact surface of class VII
whose minimal model has $b_2=0$.
Then $\Br(S)=\Br'(S)$.
\end{proposition}

Inoue showed that there are precisely three types
of Inoue surfaces, which are denoted by
$S_M$, $S^+_{N,p,q,r;t}$, and $S^-_{N,p,q,r}$.
Let us take a closer look at the first type $S_M$.
The parameter $M=(m_{ij})$ is a matrix $M\in\SL_3(\ZZ)$ that,
viewed as a complex matrix, has one real eigenvalue $\alpha>1$
and two nonreal eigenvalues $\beta,\bar{\beta}\in\CC$.
The surfaces $S=S_M$ has universal covering $\foH\times\CC$,
and the fundamental group $\pi_1(S)$ is a split extension
\begin{equation}
\label{inoue extension}
0\lra \ZZ^{\oplus 3}\lra \pi_1(S)\lra \ZZ\lra 0.
\end{equation}
The generator $g\in\ZZ$ of the quotient acts on
elements $x\in\ZZ^{\oplus 3}$ via
$gxg^{-1}=Mx$.
It requires some notation to describe
the action of  $\pi_1(S)$ on the universal covering
$\tilde{S}=\foH\times\CC$. I  do not want to write this down here,
and refer to  \cite{Inoue 1974}, Section 2. However, note  that the condition
on the eigenvalues of $M$ ensures that this action is properly discontinuous.

One easily computes that $H_1(S_M,\ZZ)=\pi_1^\ab(S_M)$ is isomorphic
to the abelian group $\ZZ\oplus\coker(M-\id)$, hence $\Br(S_M)=\Br'(S_M)=\coker(M-\id)$.
To obtain examples of matrices $M\in\SL_3(\ZZ)$ with required properties,
just choose a polynomial $p(T)=T^3+a_2T^2+a_1T+a_0$ with
integral coefficients admitting only one real root $\alpha$, and $\alpha>1$.
This holds, for example, if  $a_0\ll 0$ with $a_1,a_2$ fixed.
Now
$$
M=\begin{pmatrix}
0&0& -a_0\\
1&0&-a_1\\
0&1&-a_2
\end{pmatrix}
$$
is a matrix with characteristic polynomial
$\det(T-M)=p(T)$, hence has the desired properties.
The greatest common divisor $\delta_i$ of the $i$-minors
of the matrix $M-\id$ are $\delta_1=\delta_2=1$ and $\delta_3=p(1)=1+a_0+a_1+a_2$.
Hence the invariant factors for the
 submodule $\im(M-\id)\subset\ZZ^{\oplus 3}$ are $1,1,p(1)$.
The upshot is that $\Br(S_M)$ is cyclic of order $|1+a_0+a_1+a_2|$,
which could be arbitrarily large.

%===========================================================
\section{Main result and open questions}
\mylabel{main results}

In this section I gather our results on complex-analytic
surfaces.
In the following, let me call a compact smooth surface $S$
\emph{hypothetical} if it is of class VII, but its minimal model
is neither Hopf, Inoue, nor contains a global spherical shell.
According to the  \emph{GSS-Conjecture}  such surface should not exist.

\begin{theorem}
\mylabel{surface brauer}
Let $S$ be a smooth compact surface.
Suppose that $S$ is not  hypothetical. Then $\Br(S)=\Br'(S)$.
\end{theorem}

\proof
If $S$ is algebraic, this is Grothendieck's result
\cite{GB II},  Section 2.
For nonalgebraic surfaces we have to go through
Kodaira's classification \cite{Barth; Peters; Van de Ven 1984}, Section IV:
If the algebraic dimension is $a(S)=1$, then $S$ is elliptic and Proposition \ref{elliptic brauer} applies.
If $a(S)=0$, then $S$ is either a K3-surface, a 2-dimensional complex
torus, or  a surface of class VII.
Huybrechts and myself
settled the case of K3-surfaces  in \cite{Huybrechts; Schroeer 2003}, whereas
Elencwajg and  Narasimhan treated complex tori \cite{Elencwajg; Narasimhan 1983},
compare also Corollary \ref{lie brauer}.
We treated nonhypothetical class VII surfaces in Section \ref{Surfaces of class VII}.
\qed

\medskip
I want to finish the paper by stating some open problems:

\medskip\noindent

\medskip\noindent
(1) Suppose $S$ is a minimal surface of class VII with $b_2>0$.
Is it possible to prove that $H^3(S,\ZZ)$, or equivalently
$H_2(S,\ZZ)$ are torsion free, without referring to
the GSS-conjecture?
This would entail $\Br'(S)=0$.

\medskip\noindent
(2) Does $\Br(S)=\Br'(S)$ hold true for singular compact surfaces?

\medskip\noindent
(3) Suppose $S$ is a smooth compact surface, and
$P\ra S$ is a topological principal $\PGL_r(\CC)$-bundle. Under what conditions does there exist a holomorphic
structure on $P$?
There is a lot of work on the corresponding question for
principal $\GL_r(\CC)$-bundles. We refer to Br\^{i}nz\u{a}nescu's
book \cite{Brinzanescu 1996}.

%===========================================================


\begin{thebibliography}{ccccc}

\bibitem{Abe; Kopfermann 2001}
Y.\ Abe, K.\ Kopfermann:
Toroidal groups.
Lecture Notes in Math.\ 1759.
Springer, Berlin, 2001.

\bibitem{Barth; Peters; Van de Ven 1984}
W.~Barth, C.~Peters, A.~Van de Ven:
Compact complex surfaces.
Ergeb.\ Math.\  Grenzgebiete (3) 4,
Springer, Berlin, 1984.

\bibitem{Berkovic 1972}
V.\ Berkovi\v c:
The Brauer group of abelian varieties.
Funkcional Anal.\ i Prilo\v zen.\  6  (1972),  10--15.

\bibitem{Bogomolov 1976}
F.\ Bogomolov:
Classification of surfaces of class VII$_0$ with $b\sb{2}=0$.
Izv.\ Akad.\ Nauk SSSR Ser.\ Mat.\  40  (1976),  273--288.

\bibitem{Bogomolov; Landia 1990}
F.~Bogomolov, A.~Landia:
$2$-cocycles and Azumaya algebras under birational transformations of
algebraic schemes.
Compositio Math.\ 76 (1990),  1--5.

\bibitem{Bott; Tu 1986}
R.\ Bott, L.\ Tu:
Differential forms in algebraic topology.
Grad.\ Texts Math.\ 82. 
Berlin, Springer, 1986.

\bibitem{Brinzanescu 1996}
V.\ Br\^{i}nz\u{a}nescu:
Holomorphic vector bundles over compact complex surfaces.
Lecture Notes in Math.\ 1624.
Springer, Berlin, 1996.

\bibitem{Browder 1961}
W. Browder:
Torsion in $H$-spaces.
Ann.\ of Math.\ 74 (1961), 24--51.

\bibitem{Carlson; Toledo 1997}
J.\ Carlson, D.\ Toledo:
On fundamental groups of class VII surfaces.
Bull.\ London Math.\ Soc.\  29  (1997), 98--102.

\bibitem{Deligne 1977}
P.\ Deligne:
Cohomologie \'etale.  SGA 4$\frac{1}{2}$.
Lecture Notes in Math.\ 569.
Springer, Berlin, 1977.

\bibitem{Detloff; Grauert 1994}
G.\ Dethloff, H.\ Grauert:
Seminormal complex spaces.
In: H.\ Grauert, T.\ Peternell, R.\ Remmert eds.),
Several complex variables. VII, pages 183--220.
Encyclopaedia of Mathematical Sciences 74.
Springer, Berlin, 1994.

\bibitem{Dloussky 1984}
G.\ Dloussky:
Structure des surfaces de Kato.
M\'em.\ Soc.\ Math.\ France  14 (1984).

\bibitem{Dloussky; Oeljeklaus; Toma 2000}
G.\ Dloussky, K.\ Oeljeklaus, M.\  Toma:
Surfaces de la classe VII$\sb 0$ admettant un champ de vecteurs.
Comment.\ Math.\ Helv.\  75  (2000), 255--270.

\bibitem{Dloussky; Oeljeklaus; Toma 2003}
G.\ Dloussky, K.\ Oeljeklaus, M.\  Toma:
Class $\rm VII\sb 0$ surfaces with $b\sb 2$ curves.
Tohoku Math.\ J.\ 55  (2003), 283--309.

\bibitem{Edidin; Hassett; Kresch; Vistoli 1999}
D.~Edidin, B.~Hassett, A.~Kresch, A.~Vistoli:
Brauer groups and quotient stacks.
Amer.\ J.\ Math.\ 123 (2001), 761--777.

\bibitem{Elencwajg; Narasimhan 1983}
G.\ Elencwajg, S.\ Narasimhan:
Projective bundles on a complex torus.
J.\ Reine Angew.\ Math.\ 340 (1983), 1--5.

\bibitem{Freitag; Kiehl 1988}
E.\ Freitag, R.\ Kiehl:
\'Etale cohomology and the Weil conjecture.
Ergebnisse der Mathematik und ihrer Grenzgebiete (3)  13.
Springer, Berlin, 1988.

\bibitem{Gabber 1981}
O.\ Gabber:
Some theorems on Azumaya algebras.
In: M.~Kervaire, M.~Ojanguren (eds.), Groupe de Brauer, pp. 129--209,
Lecture Notes in Math.\ 844.
Springer, Berlin, 1981.

\bibitem{Giraud 1971}
J.~Giraud:
Cohomologie non ab\'elienne. 
Grundlehren Math.\ Wiss.\ 179. Springer, Berlin, 1971.

\bibitem{Godement 1964}
R.~Godement:
Topologie alg\'ebrique et  th\'eorie des faisceaux.
Hermann, Paris, 1964.

\bibitem{GB I}
A.\ Grothendieck:
Le groupe de Brauer I. 
In: J.~Giraud (ed.) et al.: Dix expos\'es sur la cohomologie des
sch\'emas, pp.\ 46--66.
North-Holland, Amsterdam, 1968.

\bibitem{GB II}
A.\ Grothendieck:
Le groupe de Brauer II.
In: J.~Giraud (ed.) et al.: Dix expos\'es sur la cohomologie des
sch\'emas, pp.\ 67--87.
North-Holland, Amsterdam, 1968.

\bibitem{GB III}
A.~Grothendieck:
Le groupe de Brauer III.
In: J.~Giraud (ed.) et al.: Dix expos\'es sur la cohomologie des
sch\'emas, pp.\ 88--189.
North-Holland, Amsterdam, 1968.

\bibitem{SGA 4c}
A.\ Grothendieck et al.:
Th\'eorie des topos et cohomologie \'etale. Tome 3.
Lect.\ Notes   Math.\ 305.
Springer, Berlin, 1973.

\bibitem{Hamm 1983}
H.\ Hamm:
Zum Homotopietyp Steinscher R\"aume.
J.\ Reine Angew.\ Math.\  338  (1983), 121--135.

\bibitem{Hironaka 1975}
H.\ Hironaka:
Flattening theorem in complex-analytic geometry.
Amer.\ J.\ Math.\  97  (1975), 503--547.

\bibitem{Hochschild; Serre 1953}
G.\ Hochschild, J.-P.\ Serre:
Cohomology of group extensions.
Trans.\ Amer.\ Math.\ Soc.\ 74 (1953), 110--134.

\bibitem{Hoobler 1972}
R.\ Hoobler:
Brauer groups of abelian schemes.
Ann.\ Sci.\ \'Ecole Norm.\ Sup.\  5  (1972), 45--70.

\bibitem{Huppert 1967}
B.\ Huppert:
Endliche Gruppen I.
Grundlehren Math.\ Wiss.\ 134
Springer, Berlin, 1967.

\bibitem{Huybrechts; Schroeer 2003}
D.\ Huybrechts, S.\ Schr\"oer:
The Brauer group for analytic K3-surfaces.
Int.\ Math.\ Res.\ Not.\ 50 (2003), 2687--2698.

\bibitem{Inoue 1974}
M.\ Inoue:
On surfaces of Class ${\rm VII}\sb{0}$.
Invent.\ Math.\  24  (1974), 269--310.

\bibitem{Iversen 1976}
B.\ Iversen:
Brauer group of a linear algebraic group.
J.\ Algebra  42  (1976), 295--301.

\bibitem{Iversen 1986}
B.\ Iversen:
Cohomology of sheaves. 
Springer, Berlin, 1986.

\bibitem{Iwasawa 1949}
K.\ Iwasawa:
On some types of topological groups.
Ann.\ of Math.\ 50, (1949), 507--558.

\bibitem{Karpilovsky 1987}
G.\ Karpilovsky:
The Schur multiplier.
London Math.\ Soc.\ Monogr.\  2.
Clarendon Press,  New York, 1987.

\bibitem{Kato 1975}
M.~Kato:
Topology of Hopf surfaces.
J.\ Math.\ Soc.\ Japan 27 (1975), 222--238.

\bibitem{Kato 1989}
M.~Kato:
Erratum to ``Topology of Hopf surfaces''.
J.\ Math.\ Soc.\ Japan 41 (1989),  173--174.

\bibitem{Kazama 1984}
H.\ Kazama:
$\bar \partial$ cohomology of $(H,\,C)$-groups.
Publ.\ Res.\ Inst.\ Math.\ Sci.\ 20 (1984),  297--317.

\bibitem{Kodaira 1963}
K.\ Kodaira:
On compact analytic surfaces II.
Ann.\ of Math.\  77 (1963), 563--626.

\bibitem{Kodaira 1966}
K.~Kodaira:
On the structure of compact complex analytic surfaces  II.
Am.\ J.\ Math.\ 88 (1966), 682--721.

\bibitem{Laufer 1971}
H.\ Laufer:
Normal two-dimensional singularities.
Annals of Mathematics Studies 71.
Princeton University Press, Princeton, 1971.

\bibitem{Li; Yau; Zheng 1990}
J.\ Li, S.-T.\ Yau, F.\ Zheng:
A simple proof of Bogomolov's theorem on class ${\rm VII}\sb 0$
surfaces with $b\sb 2=0$.
Illinois J.\ Math.\  34  (1990), 217--220.

\bibitem{Li; Yau; Zheng 1994}
J.\ Li, S.-T.\ Yau, F.\ Zheng:
On projectively flat Hermitian manifolds.
Comm.\ Anal.\ Geom.\ 2 (1994), 103--109.

\bibitem{Mall 1991}
D.\ Mall:
The cohomology of line bundles on Hopf manifolds.
Osaka J.\ Math.\ 28 (1991), 999--1015.

\bibitem{Maurer 1978}
J.\ Maurer:
Aufl\"osung der Entartungen holomorpher Abbildungen
zwischen zweidimensionalen Mannigfaltigkeiten.
Math.\ Ann.\  234  (1978), 89--95.

\bibitem{Milne 1980}
J.\ Milne: \'Etale cohomology. 
Princeton Mathematical Series, 33.
Princeton University Press, Princeton, 1980. 

\bibitem{Narasimhan 1967}
R.\ Narasimhan:
On the homology groups of Stein spaces.
Invent.\ Math.\ 2 (1967) 377--385.

\bibitem{Nori 1986}
M.\ Nori:
On the lattice of transcendental cycles on an elliptic surface.
Math.\ Z.\  193  (1986), 105--112.

\bibitem{Schneebeli 1979}
H.\ Schneebeli:
Group extensions whose profinite completion is exact.
Arch.\ Math.\ (Basel)  31  (1978/79), 244--253.

\bibitem{Schroeer 2001}
S.\ Schr\"oer:
There are enough Azumaya algebras on surfaces.
Math.\ Ann.\ 321 (2001), 439--454.

\bibitem{Schroeer 2003}
S.\ Schr\"oer:
The bigger Brauer group is really big.
J.\ Algebra 262 (2003), 210--225.

\bibitem{Serre 1955}
J.-P.\ Serre:
G\'eom\'etrie alg\'ebrique et g\'eom\'etrie analytique. 
Ann.\ Inst.\ Fourier  6  (1955--1956), 1--42.

\bibitem{Serre 1965}
J.-P.\ Serre:
Alg\`ebre locale. Multiplicit\'es.
Lect.\ Notes  Math.\ 11.
Springer, Berlin, 1965.

\bibitem{Serre 1972}
J.-P.\ Serre:
Cohomologie galoisienne.
Fifth edition. Lect.\ Notes Math.\ 5.
Springer, Berlin, 1994.

\bibitem{Siu 1969}
Y.-T.\ Siu:
Analytic sheaf cohomology groups of dimension n of n-dimensional complex spaces.
Trans.\ Amer.\ Math.\ Soc.\  143  (1969), 77--94.

\bibitem{Teleman 1994}
A.\ Teleman:
Projectively flat surfaces and Bogomolov's theorem
on class ${\rm VII}\sb 0$ surfaces.
Internat.\ J.\ Math.\  5  (1994), 253--264.

\bibitem{Voisin 2002}
C.\ Voisin:
Hodge theory and complex algebraic geometry. I. 
Cambridge Studies in Advanced Mathematics 76. 
Cambridge University Press, Cambridge, 2002.

\bibitem{Xiao 1991}
G.\ Xiao:
$\pi\sb 1$ of elliptic and hyperelliptic surfaces.
Internat.\ J.\ Math.\ 2 (1991),  599--615.


\end{thebibliography}
\end{document}